\documentclass[12pt]{amsart}
\usepackage{amssymb}
\begin{document}
\theoremstyle{plain}
\newtheorem{thm}{Theorem}[section]
\newtheorem{theorem}[thm]{Lemma}
\newtheorem{lemma}[thm]{Lemma}
\newtheorem{corollary}[thm]{Corollary}
\newtheorem{proposition}[thm]{Proposition}
\newtheorem{addendum}[thm]{Addendum}
\theoremstyle{definition}
\newtheorem{notations}[thm]{Notations}
\newtheorem{remark}[thm]{Remark}
\newtheorem{remarks}[thm]{Remarks}
\newtheorem{definition}[thm]{Definition}
\newtheorem{claim}[thm]{Claim}
\newtheorem{assumption}[thm]{Assumption}
\numberwithin{equation}{section}
\newcommand{\zar}{{\rm zar}}
\newcommand{\an}{{\rm an}}
\newcommand{\red}{{\rm red}}
\newcommand{\codim}{{\rm codim}}
\newcommand{\rank}{{\rm rank}}
\newcommand{\Pic}{{\rm Pic}}
\newcommand{\Div}{{\rm Div}}
\newcommand{\Hom}{{\rm Hom}}
\newcommand{\im}{{\rm im}}
\newcommand{\Spec}{{\rm Spec}}
\newcommand{\sing}{{\rm sing}}
\newcommand{\reg}{{\rm reg}}
\newcommand{\Char}{{\rm char}}
\newcommand{\Tr}{{\rm Tr}}
\newcommand{\Gal}{{\rm Gal}}
\newcommand{\Min}{{\rm Min \ }}
\newcommand{\Max}{{\rm Max \ }}
\newcommand{\soplus}[1]{\stackrel{#1}{\oplus}}
\newcommand{\dlog}{{\rm dlog}\,}    
\newcommand{\limdir}[1]{{\displaystyle{\mathop{\rm
lim}_{\buildrel\longrightarrow\over{#1}}}}\,}
\newcommand{\liminv}[1]{{\displaystyle{\mathop{\rm
lim}_{\buildrel\longleftarrow\over{#1}}}}\,}
\newcommand{\boxtensor}{{\Box\kern-9.03pt\raise1.42pt\hbox{$\times$}}}
\newcommand{\sext}{\mbox{${\mathcal E}xt\,$}}
\newcommand{\shom}{\mbox{${\mathcal H}om\,$}}
\newcommand{\coker}{{\rm coker}\,}
\renewcommand{\iff}{\mbox{ $\Longleftrightarrow$ }}
\newcommand{\onto}{\mbox{$\,\>>>\hspace{-.5cm}\to\hspace{.15cm}$}}
\catcode`\@=11
\def\opn#1#2{\def#1{\mathop{\kern0pt\fam0#2}\nolimits}}
\def\bold#1{{\bf #1}}%
\def\underrightarrow{\mathpalette\underrightarrow@}
\def\underrightarrow@#1#2{\vtop{\ialign{$##$\cr
 \hfil#1#2\hfil\cr\noalign{\nointerlineskip}%
 #1{-}\mkern-6mu\cleaders\hbox{$#1\mkern-2mu{-}\mkern-2mu$}\hfill
 \mkern-6mu{\to}\cr}}}
\let\underarrow\underrightarrow
\def\underleftarrow{\mathpalette\underleftarrow@}
\def\underleftarrow@#1#2{\vtop{\ialign{$##$\cr
 \hfil#1#2\hfil\cr\noalign{\nointerlineskip}#1{\leftarrow}\mkern-6mu
 \cleaders\hbox{$#1\mkern-2mu{-}\mkern-2mu$}\hfill
 \mkern-6mu{-}\cr}}}
\let\amp@rs@nd@\relax
\newdimen\ex@
\ex@.2326ex
\newdimen\bigaw@
\newdimen\minaw@
\minaw@16.08739\ex@
\newdimen\minCDaw@
\minCDaw@2.5pc
\newif\ifCD@
\def\minCDarrowwidth#1{\minCDaw@#1}
\newenvironment{CD}{\@CD}{\@endCD}
\def\@CD{\def\A##1A##2A{\llap{$\vcenter{\hbox
 {$\scriptstyle##1$}}$}\Big\uparrow\rlap{$\vcenter{\hbox{%
$\scriptstyle##2$}}$}&&}%
\def\V##1V##2V{\llap{$\vcenter{\hbox
 {$\scriptstyle##1$}}$}\Big\downarrow\rlap{$\vcenter{\hbox{%
$\scriptstyle##2$}}$}&&}%
\def\={&\hskip.5em\mathrel
 {\vbox{\hrule width\minCDaw@\vskip3\ex@\hrule width
 \minCDaw@}}\hskip.5em&}%
\def\verteq{\Big\Vert&&}%
\def\noarr{&&}%
\def\vspace##1{\noalign{\vskip##1\relax}}\relax\let\amp@rs@nd@&\iffalse}\fi
 \CD@true\vcenter\bgroup\relax\let\\=\cr\iffalse}\fi\tabskip\z@skip\baselineskip20\ex@
 \lineskip3\ex@\lineskiplimit3\ex@\halign\bgroup
 &\hfill$\m@th##$\hfill\cr}
\def\@endCD{\cr\egroup\egroup}
\def\>#1>#2>{\amp@rs@nd@\setbox\z@\hbox{$\scriptstyle
 \;{#1}\;\;$}\setbox\@ne\hbox{$\scriptstyle\;{#2}\;\;$}\setbox\tw@
 \hbox{$#2$}\ifCD@
 \global\bigaw@\minCDaw@\else\global\bigaw@\minaw@\fi
 \ifdim\wd\z@>\bigaw@\global\bigaw@\wd\z@\fi
 \ifdim\wd\@ne>\bigaw@\global\bigaw@\wd\@ne\fi
 \ifCD@\hskip.5em\fi
 \ifdim\wd\tw@>\z@
 \mathrel{\mathop{\hbox to\bigaw@{\rightarrowfill}}\limits^{#1}_{#2}}\else
 \mathrel{\mathop{\hbox to\bigaw@{\rightarrowfill}}\limits^{#1}}\fi
 \ifCD@\hskip.5em\fi\amp@rs@nd@}
\def\<#1<#2<{\amp@rs@nd@\setbox\z@\hbox{$\scriptstyle
 \;\;{#1}\;$}\setbox\@ne\hbox{$\scriptstyle\;\;{#2}\;$}\setbox\tw@
 \hbox{$#2$}\ifCD@
 \global\bigaw@\minCDaw@\else\global\bigaw@\minaw@\fi
 \ifdim\wd\z@>\bigaw@\global\bigaw@\wd\z@\fi
 \ifdim\wd\@ne>\bigaw@\global\bigaw@\wd\@ne\fi
 \ifCD@\hskip.5em\fi
 \ifdim\wd\tw@>\z@
 \mathrel{\mathop{\hbox to\bigaw@{\leftarrowfill}}\limits^{#1}_{#2}}\else
 \mathrel{\mathop{\hbox to\bigaw@{\leftarrowfill}}\limits^{#1}}\fi
 \ifCD@\hskip.5em\fi\amp@rs@nd@}
\newenvironment{CDS}{\@CDS}{\@endCDS}
\def\@CDS{\def\A##1A##2A{\llap{$\vcenter{\hbox
 {$\scriptstyle##1$}}$}\Big\uparrow\rlap{$\vcenter{\hbox{%
$\scriptstyle##2$}}$}&}%
\def\V##1V##2V{\llap{$\vcenter{\hbox
 {$\scriptstyle##1$}}$}\Big\downarrow\rlap{$\vcenter{\hbox{%
$\scriptstyle##2$}}$}&}%
\def\={&\hskip.5em\mathrel
 {\vbox{\hrule width\minCDaw@\vskip3\ex@\hrule width
 \minCDaw@}}\hskip.5em&}
\def\verteq{\Big\Vert&}
\def\novarr{&}
\def\noharr{&&}
\def\SE##1E##2E{\slantedarrow(0,18)(4,-3){##1}{##2}&}
\def\SW##1W##2W{\slantedarrow(24,18)(-4,-3){##1}{##2}&}
\def\NE##1E##2E{\slantedarrow(0,0)(4,3){##1}{##2}&}
\def\NW##1W##2W{\slantedarrow(24,0)(-4,3){##1}{##2}&}
\def\slantedarrow(##1)(##2)##3##4{%
\thinlines\unitlength1pt\lower 6.5pt\hbox{\begin{picture}(24,18)%
\put(##1){\vector(##2){24}}%
\put(0,8){$\scriptstyle##3$}%
\put(20,8){$\scriptstyle##4$}%
\end{picture}}}
\def\vspace##1{\noalign{\vskip##1\relax}}\relax\let\amp@rs@nd@&\iffalse}\fi
 \CD@true\vcenter\bgroup\relax\let\\=\cr\iffalse}\fi\tabskip\z@skip\baselineskip20\ex@
 \lineskip3\ex@\lineskiplimit3\ex@\halign\bgroup
 &\hfill$\m@th##$\hfill\cr}
\def\@endCDS{\cr\egroup\egroup}
\newdimen\TriCDarrw@
\newif\ifTriV@
\newenvironment{TriCDV}{\@TriCDV}{\@endTriCD}
\newenvironment{TriCDA}{\@TriCDA}{\@endTriCD}
\def\@TriCDV{\TriV@true\def\TriCDpos@{6}\@TriCD}
\def\@TriCDA{\TriV@false\def\TriCDpos@{10}\@TriCD}
\def\@TriCD#1#2#3#4#5#6{%
\setbox0\hbox{$\ifTriV@#6\else#1\fi$}
\TriCDarrw@=\wd0 \advance\TriCDarrw@ 24pt
\advance\TriCDarrw@ -1em
\def\SE##1E##2E{\slantedarrow(0,18)(2,-3){##1}{##2}&}
\def\SW##1W##2W{\slantedarrow(12,18)(-2,-3){##1}{##2}&}
\def\NE##1E##2E{\slantedarrow(0,0)(2,3){##1}{##2}&}
\def\NW##1W##2W{\slantedarrow(12,0)(-2,3){##1}{##2}&}
\def\slantedarrow(##1)(##2)##3##4{\thinlines\unitlength1pt
\lower 6.5pt\hbox{\begin{picture}(12,18)%
\put(##1){\vector(##2){12}}%
\put(-4,\TriCDpos@){$\scriptstyle##3$}%
\put(12,\TriCDpos@){$\scriptstyle##4$}%
\end{picture}}}
\def\={\mathrel {\vbox{\hrule
   width\TriCDarrw@\vskip3\ex@\hrule width
   \TriCDarrw@}}}
\def\>##1>>{\setbox\z@\hbox{$\scriptstyle
 \;{##1}\;\;$}\global\bigaw@\TriCDarrw@
 \ifdim\wd\z@>\bigaw@\global\bigaw@\wd\z@\fi
 \hskip.5em
 \mathrel{\mathop{\hbox to \TriCDarrw@
{\rightarrowfill}}\limits^{##1}}
 \hskip.5em}
\def\<##1<<{\setbox\z@\hbox{$\scriptstyle
 \;{##1}\;\;$}\global\bigaw@\TriCDarrw@
 \ifdim\wd\z@>\bigaw@\global\bigaw@\wd\z@\fi
 \mathrel{\mathop{\hbox to\bigaw@{\leftarrowfill}}\limits^{##1}}
 }
 \CD@true\vcenter\bgroup\relax\let\\=\cr\iffalse}\fi
 \tabskip\z@skip\baselineskip20\ex@
 \lineskip3\ex@\lineskiplimit3\ex@
 \ifTriV@
 \halign\bgroup
 &\hfill$\m@th##$\hfill\cr
#1&\multispan3\hfill$#2$\hfill&#3\\
&#4&#5\\
&&#6\cr\egroup%
\else
 \halign\bgroup
 &\hfill$\m@th##$\hfill\cr
&&#1\\%
&#2&#3\\
#4&\multispan3\hfill$#5$\hfill&#6\cr\egroup
\fi}
\def\@endTriCD{\egroup}
\newcommand{\sA}{{\mathcal A}}
\newcommand{\sB}{{\mathcal B}}
\newcommand{\sC}{{\mathcal C}}
\newcommand{\sD}{{\mathcal D}}
\newcommand{\sE}{{\mathcal E}}
\newcommand{\sF}{{\mathcal F}}
\newcommand{\sG}{{\mathcal G}}
\newcommand{\sH}{{\mathcal H}}
\newcommand{\sI}{{\mathcal I}}
\newcommand{\sJ}{{\mathcal J}}
\newcommand{\sK}{{\mathcal K}}
\newcommand{\sL}{{\mathcal L}}
\newcommand{\sM}{{\mathcal M}}
\newcommand{\sN}{{\mathcal N}}
\newcommand{\sO}{{\mathcal O}}
\newcommand{\sP}{{\mathcal P}}
\newcommand{\sQ}{{\mathcal Q}}
\newcommand{\sR}{{\mathcal R}}
\newcommand{\sS}{{\mathcal S}}
\newcommand{\sT}{{\mathcal T}}
\newcommand{\sU}{{\mathcal U}}
\newcommand{\sV}{{\mathcal V}}
\newcommand{\sW}{{\mathcal W}}
\newcommand{\sX}{{\mathcal X}}
\newcommand{\sY}{{\mathcal Y}}
\newcommand{\sZ}{{\mathcal Z}}
\newcommand{\A}{{\mathbb A}}
\newcommand{\B}{{\mathbb B}}
\newcommand{\C}{{\mathbb C}}
\newcommand{\D}{{\mathbb D}}
\newcommand{\E}{{\mathbb E}}
\newcommand{\F}{{\mathbb F}}
\newcommand{\G}{{\mathbb G}}
\newcommand{\HH}{{\mathbb H}}
\newcommand{\I}{{\mathbb I}}
\newcommand{\J}{{\mathbb J}}
\newcommand{\M}{{\mathbb M}}
\newcommand{\N}{{\mathbb N}}
\renewcommand{\P}{{\mathbb P}}
\newcommand{\Q}{{\mathbb Q}}
\newcommand{\R}{{\mathbb R}}
\newcommand{\T}{{\mathbb T}}
\newcommand{\U}{{\mathbb U}}
\newcommand{\V}{{\mathbb V}}
\newcommand{\W}{{\mathbb W}}
\newcommand{\X}{{\mathbb X}}
\newcommand{\Y}{{\mathbb Y}}
\newcommand{\Z}{{\mathbb Z}}
\title[On the Shafarevich conjecture]{On the Shafarevich conjecture for surfaces of general
type over function fields}
\author{Egor Bedulev}
\address{Independent Moscow University, Moscow, Russia}
\email{bedulev@mccme.ru}
\thanks{This work has been partly supported by the DFG Forschergruppe ``Arithmetik
und Geometrie''} 
\author{Eckart Viehweg}
\address{Universit\"at GH Essen, FB6 Mathematik, 45117 Essen, Germany}
\email{ viehweg@uni-essen.de}
\maketitle

Let $B$ be a projective algebraic curve over an algebraically
closed field $k$ of characteristic zero, and let $S \subseteq B$
be a finite set of points. Shafarevich's conjecture for families
of curves, proved by Parshin and Arakelov (see \cite{A}), states
that
\begin{quote}
(I) There are only finitely many isomorphism classes of smooth
non-isotrivial families of curves over $B-S$.
\end{quote}
\begin{quote}
(II) If $2q - 2 + \# S \leq 0$, then there are no such families.
\end{quote}
Similar questions can be asked for smooth families of higher
dimensional manifolds. (II) was recently verified for families
of minimal surfaces of general type and for canonically
polarized manifolds, by Kov\'acs \cite{K 2}, \cite{K 1}, \cite{K 3},
Migliorini \cite{M 1} and Zhang \cite{Z 1}. As a byproduct we
reprove their result, but the reader familiar with their
articles will recognize strong similarities between their and
our approach.

Unfortunately our method does not work for families
whose relative dualizing sheaf is big and numerically effective
on the smooth fibers. In this case, for non-isotrivial families over
$\P^1$ or over an elliptic curve, Kov\'acs (Thesis, Utah, 1995,
see \cite{K 2}) showed the existence of at least one degenerate
fibre.

In the higher dimensional case (I) might be too much to hope
for. For fixed $B$ and $S$ there exist non-trivial deformations
of families of abelian varieties over $B-S$ (see \cite{F}).
So (I) splits up in two questions: boundedness and rigidity.

To be more precise, let us fix some polynomial $h \in \Q [t]$
with $h (\Z) \subset \Z$. If ${\rm deg} (h) =2$, we define $M_h$
to be the moduli scheme of minimal surfaces $S$ of general type
with $h(\mu) = \chi (\omega^{\mu}_{S})$, or allowing
singularities, the moduli scheme of canonically polarized normal
surfaces with at most rational double points and with Hilbert
polynomial $h$ (see \cite{V 1}, for example).

If ${\rm deg} (h) > 2$, we denote by $M_h$ the moduli scheme of
canonically polarized manifolds, with Hilbert polynomial $h$.
In both cases (I) should be replaced by two sub-problems:
\begin{quote}
(B) The non-trivial morphisms $B - S \to M_h$, which are induced
by smooth projective maps $g_0 : Y_0 \to B - S$, are
parameterized by some scheme of finite type.
\end{quote}
\begin{quote}
(R) Under which additional conditions are the morphisms $B - S
\to M_h$ in (B) rigid.
\end{quote}

The corresponding questions, for abelian varieties and their
moduli scheme, were solved by Faltings \cite{F}. In this note we
prove (B) for surfaces of general type, and for canonically
polarized manifolds, in case $S = \emptyset$. The only obstruction to
extend (B) to arbitrary families of canonically polarized
manifolds, is the lack of a proof for the existence of relative
minimal models for semi-stable families of such varieties over
curves. We have nothing to contribute to problem (R).

For smooth families $g : Y \to B$ of surfaces, i.e.
for $S = \emptyset$, (B) has been proved by the first author.

This note benefited from discussions between the first author
and A. Parshin and between the second author and Qi Zhang. We
thank both of them, and K. Oguiso for remarks and comments.

\section{Families of canonically polarized manifolds and of
surfaces of general type}

Let $B$ be a curve and $Y$ a variety of dimension $n+1$, both
non-singular, projective and defined over an algebraically
closed field $k$ of characteristic zero. Let $g : Y \to B$ be a
flat morphism with connected fibers $Y_b = g^{-1} (b)$. Let $S
\subset B$ be a finite set of points and $D = g^{-1} (S)$.
Frequently we will denote the divisors $\sum_{s \in S} s$ and
$\sum_{s \in S} g^{-1} (s)$ by $S$ and $D$, as well.

We want to study both, families of surfaces of general type and
of canonically polarized manifolds, and $D$ will be supposed to
be the set of bad fibers. To cover both cases we formulate the
following assumption:

\begin{assumption} \label{fam-ass} \
\begin{enumerate}
\item[a)] $g|_{Y-D} : Y - D \>>> B - S$ is smooth and
$\omega_{Y/B}|_{Y-D}$ is relatively semi-ample, i.e. for some $\nu
\gg 0$ the natural map
$$
\Phi_{\nu} : g^* g_* \omega^{\nu}_{Y/B} \>>> \omega^{\nu}_{Y/B}
$$
is surjective over $Y-D$.
\item[b)] The fibers $Y_b$ should be manifolds of general type
and the $\nu$-canonical map should contract at most curves.
Hence if
$$
\pi_{\nu} : Y - D \twoheadrightarrow  V \subseteq \P (g_* \omega^{\nu}_{Y/B}
|_{B-S} )
$$
is the morphism, induced by $\Phi_{\nu}$, then $\pi_{\nu}$
should be birational and the maximal fibre dimension of
$\pi_{\nu}$ should be one, for some $\nu \gg 0$, satisfying the
condition a).
\end{enumerate}
\end{assumption}

Recall from \cite{EV 2}, \S \ 2, or \cite{V 1}, \S \ 5, that for
an invertible sheaf $\sL $ on a normal projective variety $F$
the integer $e (\sL) $ is defined to be the smallest positive
integer $e$ for which the multiplier sheaf $\omega_F \left\{ -
\frac{D}{e} \right\}$ is isomorphic to $\omega_F$, for all
divisors $D$ of global sections of $\sL$. For $e (\sL)$ to
exist, one has to require $F$ to have at most rational
singularities. Then, by \cite{EV 2}, 2.3, or \cite{V 1}, 5.11,
5.12 and 5.21, $e (\sL)$ exists and in some cases there are
explicit bounds:

\begin{lemma} \label{multiplier}
Assume $F$ is projective with at most rational singularities and
let $\sL$ be an ample invertible sheaf on $F$.
\begin{enumerate}
\item[a)] If $F$ is non-singular and $\sL$ very ample, then
$e(\sL) \leq c_1 (\sL)^{\dim F} +1$.
\item[b)] In general, let $\sL$ be very ample and assume that
there exists a desingularization $\sigma : F' \to F$ and an
effective divisor $E$ such that $\sigma^*  \sL \otimes \sO_{F'}
(-E)$ is very ample. Then $e (\sL) \leq c_1 (\sL)^{\dim F} +1$.
\item[c)] If $F$ has rational Gorenstein singularities and if $Z
= F \times \ldots \times F$ and $\sM = \otimes pr^{*}_{i} \sL$,
then $e (\sM) = e (\sL). $
\end{enumerate}
\end{lemma}

Recall that $g: Y \to B$ is called isotrivial, if for some
variety $F$, defined over $k$, there is a birational map
$Y{\times_B} \overline{k(B)} \to F{\times_k} \overline{k(B)}$.

In \cite{V 1}, 2.7 and 2.9, we defined a locally free sheaf $\sG
$ on $B$ to be numerically effective (or nef), if for all $\mu
>0$ and for a point $p \in B$ the sheaf
$$
S^{\mu} (\sG) \otimes \sO_B (p)
$$
is ample. We will use:

\begin{proposition} \label{pos}\
\begin{enumerate}
\item[a)] If $g: Y \to B$ is a projective morphism between
non-singular varieties, then $g_* \omega^{\nu}_{Y/B}$ is nef,
for all $\nu >0$.
\item[b)] If $g$ satisfies the assumptions made in
\ref{fam-ass}, a), if the general fibre of $g$ is of general type
and if $g$ is non-isotrivial, then $\kappa (\omega_{Y/B}) = n+1$
and ${\rm det} (g_* \omega^{\eta}_{Y/B})$ is ample, provided
$g_* \omega^{\eta}_{Y/B} \neq 0$ and $\eta > 1$.
\end{enumerate}
\end{proposition}
a) is a special case of \cite{V 2}, Theorem III, and b) can be
found in \cite{V 3}, Theorem II. In fact, there is the ampleness of
${\rm det} (g_* \omega^{\eta}_{Y/B})$ is shown for some $\eta$,
but as in \cite{EV 2}, 3.1, or \cite{V 1} this implies that $g_*
\omega^{\eta}_{Y/B}$ is ample for all $\eta >1$.

\begin{thm} \label{bounds}
Let $g : Y \to B$ be a non-isotrivial morphism, satisfying the
assumptions made in \ref{fam-ass}. Let us fix some $\nu >1$,
such that \ref{fam-ass} a) and b) hold true. We write
\begin{itemize}
\item $q = g (B)$ for the genus of $B$.
\item $s = \# S$ for the number of degenerate fibers.
\item $e (\nu) = e (\omega^{\nu}_{F})$ for some general fibre
$F$ of $g$.
\item $r (\nu) = {\rm rank} (g_* \omega^{\nu}_{Y/B})$.
\end{itemize}
Then one has:
\begin{enumerate}
\item[a)] $2q - 2 + s > 0$
\item[b)] If $g : Y \to B$ is semi-stable, then $$n \cdot (2q - 2
+s) \cdot \nu \cdot e (\nu) \cdot r (\nu) \geq {\rm deg} (g_*
\omega^{\nu}_{Y/B}).$$
\item[c)] In general $$ (n \cdot (2q - 2 +s) +s) \cdot \nu \cdot
e (\nu) \cdot r (\nu) \geq {\rm deg} (g_* \omega^{\nu}_{Y/B}).$$
\end{enumerate}
\end{thm}
a) will follow from \ref{pos} and b). This part of theorem
\ref{bounds} is due to  Kov\'acs \cite{K 2}, \cite{K 1}, \cite{K
3} and to Migliorini \cite{M 1}. Qi Zhang \cite{Z 1} gave an
elegant proof in case $B$ is an elliptic curve. His proof easily
extends to $B = \P^1$, as he and the second author found out
discussing his result. In this note we take up their approach,
together with positivity properties of direct image sheaves, as
stated in \cite{EV 2}.

\begin{remark} \label{bounds 2}
The constants $\nu$ and $r (\nu)$ are determined, using Matsusaka's
big theorem, by the Hilbert polynomial of the fibers $Y_b$ for
$b \in B - S$. For canonically polarized manifolds, $e(\nu) \leq
\nu^n \cdot c_1 (\omega_F)^n +1$, as we have seen in
\ref{multiplier}, a).

For surfaces of general type, the
canonical model $\tilde{F}$ of $F$ has $A - D - E$
singularities, and the number of $(-2)$-curves on $F$ is
bounded by $\dim H^1 (F, \Omega^{1}_{F})$, a number
which is constant in families. One can use \ref{multiplier}, b)
to bound $e (\omega^{\nu}_{F}) =e (\omega^{\nu}_{\tilde{F}} )$.

Such a bound exists for different reasons. By \cite{V
1}, the moduli space of normal canonically polarized surfaces
with at most rational double points and
with given Hilbert polynomial is quasi-projective, in particular
it has only finitely many irreducible components. Moreover,
Koll\'ar \cite{Ko} (see also \cite{V 1}, 9.25) constructed a
finite covering $Z$ of the moduli scheme, together with a
``universal family''. By \cite{V 1}, 5.17,
$e(\omega^{\nu}_{\tilde{F}})$ is bounded on each irreducible
component of $Z$. Hence there exists some constant $e$ depending
on $h$ with $e(\omega^{\nu}_{\tilde{F}})\leq e$, for all
surfaces $F$ with Hilbert polynomial $h$.
\end{remark}

\section{Ampleness and vanishing theorems}

Let $X$ be a projective manifold, $U \subset X$ an open dense
submanifold, and let $\sL$ be an invertible sheaf on $X$.

\begin{definition} \label{def-sa}\
\begin{enumerate}
\item[a)] $\sL$ is semi-ample with respect to $U$ if for some
$\eta > 0$
$$
\iota_{\eta} : H^0  (X, \sL^{\eta}) \otimes_k \sO_X \>>>
\sL^{\eta}
$$
is surjective over $U$.
\item[b)] $\sL$ is $\ell$-ample and big with respect to $U$, if one
can choose $\eta$ in a) such that the morphism
$$
\Phi_{\eta} : U \twoheadrightarrow V \subset \P (H^0 (X,
\sL^{\eta}))
$$
induced by $\iota_{\eta}$ is proper, birational and
$$
{\rm Max} \{ \dim \Phi^{-1}_{\eta} (v) ; v \in V \} \leq l.
$$
\end{enumerate}
\end{definition}
A slight modification of the argument used to prove 6.4 in
\cite{EV 1}, yields a generalization of the
Akizuki-Kodaira-Nakano vanishing theorem, similar to the one
used in the proof of \cite{K 1}, 1.1 and \cite{K 3}, 1.1.

\begin{proposition} \label{vanishing}
Let $\sL$ be $\ell$-ample and big with respect to $U$.
Assume that for $\eta$ as in \ref{def-sa}, b), the image of $V$
of $\Phi_{\eta}$ allows a projective flat morphism
$\gamma : V \to W$ to a non-singular affine variety $W$.

Then there exists a blowing up $\tau : X' \to X$ with centers in
$\Delta = X = U$, such that $\Delta' = \tau^{-1} (X - U)$ is a
normal crossing divisor, and an effective divisor $\Gamma'$ with
$\Gamma'_{{\rm red}} \leq \Delta'$, such that for all
numerically effective invertible sheaves $\sN$ and for $p +q <
\dim X - {\rm Max} \{0,\ell-1\}$,
$$
H^p (X', \Omega^{q}_{X'} ({\rm log} \ \Delta' ) \otimes \tau^*
(\sL^{-1} \otimes \sN^{-1}) \otimes \sO_{X'} (\Gamma')) =0.
$$
\end{proposition}

\begin{proof}
If $\Delta = \emptyset$, this is \cite{EV 1}, 6.6. Hence we will
assume that $\Delta \neq 0$, and fix some $\eta$ for which the
assumption \ref{def-sa} b) on $\Phi_{\eta}$ hold true. Let $\tau
: X' \to X$ be a blowing up, such that $X'$ is non-singular,
$\Delta' = \tau^{-1} (\Delta)$ a normal crossing divisor and
such that, for $\sL' = \tau^* \sL$, the image of
$$
\iota'_{\eta} : H^0 (X', \sL^{'\eta} ) \otimes_k \sO_{X'} \>>>
\sL^{'\eta}
$$
is an invertible sheaf, isomorphic to $\sL^{'\eta}$ over $U' =
\tau^{-1} (U)$.

So ${\rm Im} (\iota'_{\eta}) = \sL^{'\eta} (- \Gamma_1)$, for
some divisor $\Gamma_1$, supported in $\Delta'$. Let
$$
\Phi' : X ' \twoheadrightarrow Z \subset \P (H^0 (X',
\sL^{'\eta} (- \Gamma_1 )))
$$
be the induced morphism. $\Phi' |_{U'}$ is a proper morphism
with image $V$. Let $I$ be the ideal sheaf of $Z - V$ and
$$J =\Phi^{'*} I/_{\rm torsion}.$$
Blowing up again, we may assume that $J =
\sO (- \Gamma_2)$ for some divisor $\Gamma_2$ with
$(\Gamma_2)_{{\rm red}} = \Delta'_{{\rm red}}$. For all $\mu \gg
0$ the sheaf $\sL^{'\mu \cdot \eta} ( - \mu \cdot \Gamma_1 - \Gamma_2)$
will be generated by global sections. For some $\mu \gg 0$, the
number $\mu \cdot \eta$ will not divide the multiplicity of any
of the components of $\mu \cdot \Gamma_1 + \Gamma_2$. Allowing
$\Delta'$ to have multiplicities, we thereby got to the
following situation:

\begin{assumption} \label{ass-van}
For some $\eta$ and some normal crossing divisor $\Delta'$ with
$$\Delta'_{{\rm red}} = X' - U' = (\Delta' - \eta \cdot \Bigl[
\frac{\Delta'}{\eta} \Bigr])_{{\rm red}},$$
the sheaf $\sL^{'\eta} (-\Delta') = \sL^{'\eta} \otimes \sO_{X'} (-
\Delta') $ is generated by global sections, the induced morphism
$\Phi' : X' \twoheadrightarrow Z$ is birational and the
fiber-dimension of $U' = \Phi^{'-1} (V) \to V$ is at most $\ell$.
Moreover there exists a projective flat morphism $\gamma : V \to W$
to a non-singular affine variety $W$. Let $\sN'$ be any
numerically effective invertible sheaf on $X'$ and choose
$[\frac{\Delta'}{\eta}] = \Gamma '$.
\end{assumption}
We will prove, by induction on $\dim X' - \dim W$, that
the assumptions \ref{ass-van} imply
\begin{gather} \label{eq-van}
H^p (X', \Omega^{q}_{X'} ({\rm log} \ \Delta') \otimes \sL^{'-1}
\otimes \sN^{'-1} \otimes \sO_{X'} (\Gamma')) =0 \\
\mbox{for \ \ } \
p + q < \dim X - {\rm Max} \{ 0,r-1\}. \notag
\end{gather}
If $\dim X' = \dim V > \dim W$
we may replace $\eta$ by $\mu \cdot \eta$ and $\Delta'$ by $\mu
\cdot \Delta'$ for $\mu \gg 0$. Thereby we are allowed to assume
that $(\sL' \otimes \sN')^{\eta} \otimes \sO_{X'} (- \Delta')$
is generated by global sections, and that for the zero divisor
$H$ of a general section of this sheaf $\Phi' (H) \cup V \to W$
is again projective and flat. Moreover $\Delta' + H$
is a normal crossing divisor and $\Gamma' |_{H} =
[\frac{\Delta'|_{H}}{\eta}]$. Therefore $H$, $\sL'|_H$ and
$\Delta'|_H$ satisfy again the assumptions \ref{ass-van}. If
$\dim V = \dim W$, we choose $H = 0$. In both cases, the morphism
$$V- \Phi'(H) \>>> W$$
is affine and $V - \Phi' (H)$ is an affine variety. $(\sL'
\otimes \sN')^{\eta}$ has a section with zero
divisor $\Delta' +H$. By \cite{EV 1}, \S \ 3,
$$
(\sL' \otimes \sN')^{-1} \otimes \sO_{X'} (\Gamma')
=(\sL' \otimes \sN')^{-1} \otimes \sO_{X'} (\Bigl[
\frac{\Delta'+H}{\eta} \Bigr])
$$
has an integrable connection which satisfies the
$E_1$-degeneration. The assumption
$\Delta'_{{\rm red}} = (\Delta' - \eta \cdot [
\frac{\Delta'}{\eta} ])_{{\rm red}}$ allows to apply \cite{EV 1},
4.12, and to obtain
$$
H^p(X', \Omega^{q}_{X'} ({\rm log} \ \Delta' + H') \otimes (\sL'
\otimes \sN')^{-1} \otimes \sO_{X'} (\Gamma')) =0,
$$
for $p+q < n - {\rm Max} \{ 0,\ell-1\}$. If $\dim V = \dim W$, we
are done. If $\dim V > \dim W$ one uses the exact sequence
$$
0 \to \Omega^{p}_{X'} ({\rm log} \ \Delta') \to \Omega^{p}_{X'}
({\rm log} \ \Delta' + H') \to \Omega^{p-1}_{H'} ({\rm log} \
\Delta' |_{H'} ) \to 0
$$
to obtain \ref{eq-van} by induction.
\end{proof}

\section{The proof of theorem \ref{bounds}}

We will assume that $D$ is a normal crossing divisor. Enlarging
$S$ and $D$ we may assume that $2q - 2 + s \geq 0$, hence
that $\omega_B (S)$ is nef. If $2q - 2 +
s =0$, then either $B$ is elliptic and $g$ smooth, or $B =
\P^1$ and $S = \{ b_1, b_2\}$. In the second case, there
exists a finite covering $\P^1 \to \P^1$, totally ramified in
$S$, such that the pullback family has stable reduction.
Altogether, part a) of \ref{bounds} follows from \ref{pos} and
part b) of \ref{bounds}.

If the fibers of $Y -D \to B - S$ are canonically polarized manifolds, the
next proposition follows from \cite{EV 2}, 2.4. Since we will
use it for families of minimal models of surfaces of general
type, as well, we will recall the proof.

\begin{proposition} \label{ample}
Under the assumptions made in \ref{bounds} let $\sN$ be an
invertible sheaf on $B$ with
${\rm deg} \ \sN < {\rm deg} \ g_* \omega^{\nu}_{Y/B}.$
Then the sheaf
$$
\omega^{\nu \cdot e (\nu) \cdot r(\nu)}_{Y/B} \otimes g^*
\sN^{-1}
$$
is 1-ample and big with respect to $Y-D$.
\end{proposition}

\begin{proof} By the assumptions a) and b) in \ref{fam-ass} it
is sufficient to show, that for some $\mu >0$
$$
S^{\mu \cdot r (\nu) \cdot e(\nu)} (g_* \omega^{\nu}_{Y/B} )
\otimes \sN^{- \mu}
$$
is generated by global sections, or that
$$
S^{r (\nu) \cdot e(\nu)} (g_* \omega^{\nu}_{Y/B}) \otimes
\sN^{-1} $$
is ample. By definition of nef, this hold true, if the sheaf
$$
S^{r(\nu) \cdot e(\nu)} (g_* \omega^{\nu}_{Y/B} ) \otimes {\rm
det} (g_* \omega^{\nu}_{Y/B} )^{-1}
$$ is nef. To this aim (see \cite{V 1}, 2.8) we can replace $B$
by a covering, unramified in $S$, and $Y$ by the pullback
family. Thereby we may assume that
$$
{\rm det} (g_* \omega^{\nu}_{Y/B} ) = \sA^{e(\nu)}
$$
for an invertible sheaf $\sA$ on $B$. For $r = r (\nu)$ let
$$
f : X = Y \times_B \ldots \times_B Y \>>> B
$$
be the $r$-fold fibre product, let $\sigma: X' \to X$ be a
desingularization, and
$$f' = f \circ \sigma : X' \to B.$$
We write $\sM = \sigma^* (\otimes^{r}_{i=1} pr^{*}_{i}
\omega_{Y/B})$. The morphism $f$ is Gorenstein and the general
fibre is non-singular. Hence there are natural injective maps
\begin{gather}
\otimes^r g_* \omega^{\nu}_{Y/B} = f_{*} \otimes^{r}_{i=1}
pr^{*}_{i} \omega^{\nu}_{Y/B} \>>> f'_{*} \sM^{\nu}
\label{eq1} \\
\mbox{and} \ \ \ f'_* \sM^{\nu -1} \otimes \omega_{X'/B} \>>> f_*
\omega^{\nu}_{X/B} = \otimes^r g_* \omega^{\nu}_{Y/B},
\label{eq2}
\end{gather}
both isomorphisms on some open dense subset of $B$. (\ref{eq1}) induces
$$
\sA^{e(\nu)} = {\rm det} (g_* \omega^{\nu}_{Y/B}) \>>> \otimes^r
g_* \omega^{\nu}_{Y/B} \>>> f'_* \sM^{\nu} ,
$$
hence a section of $\sM^{\nu} \otimes f^{'*} \sA^{-e(\nu)}$
with zero-divisor $\Gamma$. Blowing up $X'$, with centers in a
finite number of fibers, we may assume that
the image $\sM^{\nu} \otimes J$ of
$$
f^{'*} \otimes^r g_* \omega^{\nu}_{Y/B} \>>> f^{'*} f'_*
\sM^{\nu} \>>> \sM^{\nu}
$$
is invertible. The ideal sheaf $J$ is trivial in a neighborhood
of the general fibre. By \ref{pos} $g_* \omega^{\nu}_{Y/B}$ is
nef, hence $\sM^{\nu} \otimes J$ is nef, as well.

Let us write $\sL = \sM^{\nu-1} \otimes f^{'*} \sA^{-1}$. Then
$$
\sL^{e(\nu) \cdot \nu}(-\nu\cdot\Gamma)
= \sM^{e(\nu) \cdot \nu \cdot (\nu-1)}\otimes f^{'*}
\sA^{-\nu\cdot e(\nu)} \otimes \sO_{X'}(-\nu\cdot\Gamma)
= \sM^{e(\nu) \cdot \nu \cdot (\nu-1) - \nu^2}
$$
contains a nef subsheaf, isomorphic to $\sL^{e(\nu) \cdot
\nu}(-\nu\cdot\Gamma)$ in a neighborhood of the
general fibre. Moreover, by \ref{pos}, b),
$$\kappa (\sL^{e (\nu) \cdot
\nu} (-\nu \cdot \Gamma)) = \dim X'.$$
\cite{EV 2}, 1.7 implies that the sheaf $f'_* \sL \otimes
\omega_{X'/B} \{ - \frac{\Gamma}{e(\nu)} \}$ is nef. By the
definition of $e(\nu)$ and by \ref{multiplier}, c), the natural inclusion
$$
f'_* \sL \otimes \omega_{X'/B} \Bigl\{ - \frac{\Gamma}{e (\nu)} \Bigr\}
\>>> f'_* \sL \otimes \omega_{X'/B} = f'_* (\sM^{\nu-1} \otimes
\omega_{X'/B}) \otimes \sA^{-1}
$$
is an isomorphism on some open dense subset of $B$. Using
(\ref{eq2}) one obtains a nef subsheaf of
$$
(\otimes^r g_* \omega^{\nu}_{Y/B} ) \otimes \sA^{-1}
$$
of full rank, hence the latter is nef, as well. Then
$$
S^{e(\nu)} (\otimes^r g_* \omega^{\nu}_{Y/B}) \otimes
\sA^{-e(\nu)} = S^{e(\nu)} (\otimes^r g_* \omega^{\nu}_{Y/B} )
\otimes {\rm det} (g_* \omega^{\nu}_{Y/B})^{-1}
$$
as well as its quotient $S^{e (\nu) \cdot r} (g_*
\omega^{\nu}_{Y/B}) \otimes {\rm det} (g_* \omega^{\nu}_{Y/B} )^{-1}$
are nef.
\end{proof}

Let us return to the proof of \ref{bounds} b) and c). If $g$ is
semi-stable, i.e. if $D$ is reduced, we choose $\delta=0$, otherwise
$\delta=1$. Let us assume that \ref{bounds} b) or c) are wrong. Hence
\begin{equation} \label{equality}
(n\cdot (2q -2+s) +\delta \cdot s) \cdot \nu \cdot e (\nu) \cdot
r(\nu) < {\rm deg} (g_* \omega^{\nu}_{X/B})
\end{equation}
and \ref{ample} implies that for
$$
\sA = \omega_{B} (S)^n \otimes \sO_B (\delta \cdot S)
$$
the sheaf $\omega^{\nu \cdot e (\nu) \cdot r(\nu)}_{Y/B} \otimes
\sA^{-\nu \cdot e(\nu) \cdot r (\nu)}$, hence $\sL =
\omega_{Y/B} \otimes g^* \sA^{-1}$, is 1-ample and big with
respect to $Y-D$. Since we assumed $\omega_{B} (S)$ to be nef,
\ref{vanishing} implies that:

\begin{claim} \label{van-am}
There exists a blowing up $\tau : X \to Y$ with centers in $D$,
such that $\Delta = \tau^* D$ is a normal crossing divisor, and
there exists an effective divisor $\Gamma$, supported in
$\Delta$, with
$$
H^p (X, \Omega^{q}_{X} ({\rm log} \ \Delta) \otimes \tau^*
\omega^{-1}_{Y/B} \otimes g^* (\omega_B (S)^{n-m} \otimes
\sO_{B} (\delta \cdot S)) \otimes \sO_X (\Gamma)) = 0
$$
for $p + q < \dim X = n+1$ and for all $m \geq 0$.
\end{claim}
The morphism $f = g \circ \tau$ is smooth outside of $\Delta$,
and one has an exact sequence of locally free sheaves
\begin{equation} \label{seq1}
0 \to f^* \omega_B (S) \to \Omega^{1}_{X} ({\rm log} \Delta)
\to \Omega^{1}_{X/B} \to 0.
\end{equation}
Comparing the determinants one finds
$$\Omega^{n}_{X/B} = {\rm det} (\Omega^{1}_{X/B}) = \omega_X
(\Delta_{{\rm red}}) \otimes f^* \omega_B (S)^{-1} =
\omega_{X/B} (\Delta_{{\rm red}} - \Delta).$$

\begin{claim} \label{effective}
$$
H^0 (X, \Omega^{n}_{X/B} \otimes \tau^* \omega^{-1}_{Y/B}
\otimes g^* \sO_B (\delta \cdot S)) \neq 0
$$
\end{claim}

\begin{proof} If $\delta =1$, this holds true with $\tau^*
\omega^{-1}_{Y/B}$ replaced by the smaller sheaf
$\omega^{-1}_{X/B}$, since
$$
\Omega^{n}_{X/B} \otimes \omega^{-1}_{X/B} \otimes g^* \sO_B (S)
= \sO (\Delta_{{\rm red}}).
$$
If $\delta =0$, i.e. if $g$ is semi-stable, then $\Omega^{n}_{Y/B} =
\omega_{Y/B}$. By \cite{EV 1}, 3.21, one has an inclusion
$\tau^* \Omega^{1}_{Y} ({\rm log} \ D) \to \Omega^{1}_{X} ({\rm
log} \ \Delta)$, hence $\tau^* \omega_{Y/B} = \tau^*
\Omega^{n}_{Y/B} \subset \Omega^{n}_{X/B}$.
\end{proof}
The $m$-th wedge product, applied to the sequence (\ref{seq1})
induces an exact sequence
\begin{equation}
0 \to f^* \omega_B (S) \otimes \Omega^{m-1}_{X/B} \to
\Omega^{m}_{X} ({\rm log} \ \Delta) \to \Omega^{m}_{X/B} \to 0 .
\tag{$\Sigma^{\bullet}_{m}$}
\end{equation}
\ref{van-am} implies, that the $(n - m)$-th cohomology of the
sheaf in the middle of the sequence
$$
\Sigma^{\bullet}_{m} \otimes \tau^* \omega^{-1}_{Y/B} \otimes
g^* \omega^{n-m}_{B} ((n-m + \delta) \cdot S) \otimes \sO_X
(\Gamma)
$$
is zero, hence
$$
H^{n-m} (X, \Omega^{m}_{X/B} \otimes \tau^* \omega^{-1}_{Y/B}
\otimes g^* \omega^{n-m}_{B} ((n-m + \delta) \cdot S) \otimes
\sO_X (\Gamma))
$$
injects into
$$
H^{n-m+1} (X, \Omega^{m-1}_{X/B} \otimes \tau^*
\omega^{-1}_{Y/B} \otimes g^* \omega^{n-m+1}_{B} ((n-m+1+\delta)
\cdot S) \otimes \sO_X (\Gamma)).
$$
Altogether
$$
H^0 = H^0 (X, \Omega^{n}_{X/B} \otimes \tau^* \omega^{-1}_{Y/B}
\otimes g^* \sO_B (\delta \cdot S) \otimes \sO_X (\Gamma))
$$
is a subspace of
$$
H^n = H^n (X, \tau^* \omega^{-1}_{Y/B} \otimes g^* \omega^{n}_{B}
((n+\delta) \cdot S) \otimes \sO_X (\Gamma)).
$$
Using \ref{van-am} again, one finds $H^n$ and thereby $H^0$ to be zero.
Since $\Gamma$ is effective, this contradicts \ref{effective}.
Hence the inequality stated in (\ref{equality}) does not hold
true, and we obtain \ref{bounds} b) and c).

\section{Moduli schemes and boundedness}

Let $\sM_h$ denote the moduli functor of minimal
surfaces of general type, if ${\rm deg} \ h = 2$, and of
canonically polarized manifolds, if ${\rm deg} \ h > 2$, both
times with fixed Hilbert polynomial $h\in \Q[t]$. The
corresponding moduli scheme is denoted by $M_h$.

If ${\rm deg} \ h = 2$, i.e. in the surface case, Koll\'ar and
Shepherd-Barron \cite{KS} defined stable surfaces, Alexeev
\cite{Al} proved that the index of the singularities is bounded
in terms of the coefficients of the Hilbert polynomial, which
implies by \cite{Ko} that $M_h$ has a compactification
$\bar{M}_h$, parameterizing families of stable surfaces (see also
\cite{V 1}, section 9.6).

For ${\rm deg} \ h \geq 2$, those
results were generalized by Karu \cite{KK}, assuming the minimal
model conjecture for semi-stable families of $n$-folds over
curves $(MMP (n+1))$. Let again $\bar{M}_h$ be the
compactification and $\bar{\sM}_h$ the corresponding moduli
functor. The existence of $\bar{M}_h$ implies that for some
$N_0$, depending on $h$, the reflexive hull $\omega^{[N_0]}_{X}$
of $\omega^{N_0}_{X}$ is invertible for all $x \in \bar{\sM}_h
(k)$.

In both cases, for some $\eta\gg 0$ the sheaf
$\omega^{[\eta]}_{X}$ is very ample $({\rm deg} \ h > 2)$ or
semi-ample and the induced morphism the contradiction of $(-2)$
curves $({\rm deg} \ h = 2)$.

Koll\'ar \cite{Ko} (see also \cite{V 1}) has shown, that
for $\eta\gg 0$ and for some $p >0$
there exists a very ample invertible sheaf $\lambda$ on $\bar{M}_h$
with:
\begin{quote}
For $f: X \to Z \in \bar{\sM}_h (Z)$ and for the
induced morphism $\varphi : Z \to \bar{M}_h$,
$$\varphi^* \lambda
= {\rm det} (f_* \omega^{[\eta]}_{X/Z})^p.$$
\end{quote}
\begin{corollary} \label{b_s}
Assume ${\rm deg} \ h = n = 2$. Let $B$ be a projective
non-singular curve, $S \subset B$ a finite subset, and let
$g_0 : Y_0 \to B - S$ be a smooth projective morphism, whose
fibers are surfaces of general type with Hilbert polynomial $h$.
Let $\Phi : B \to \bar{M}_h$ be the induced morphism. Then ${\rm
deg} \ \Phi^* \lambda$ is bounded above by a constant, depending
on $h, g (B)$ and $\# S$.
\end{corollary}

\begin{addendum} \label{b_cp}
Assuming $MMP ({\rm deg} \ h + 1 )$, \ref{b_s} remains true for
${\rm deg} \ h > 2$ and for $g_0 : Y_0 \to B - S$ a family of
canonically polarized manifolds with Hilbert polynomial $h$.
\end{addendum}

\begin{proof} Choose a non-singular projective compactification
$Y$ of $Y_0$ such that $g_0$ extends to $g: Y \to B$. The
assumptions \ref{fam-ass} hold true for $\nu = \eta$. By
\ref{bounds} the degree of ${\rm det} (g_* \omega^{\eta}_{Y/B})$
is smaller than a constant depending on $h, g (B)$ and $\# S$.
There exists a finite covering $\gamma : C \to B$ and $f: X \to
C \in \bar{\sM}_h (C)$ which induces
$$
C \> \gamma >> B \> \Phi >> \bar{M}_h .
$$
We may assume, in addition, that $f$ has a semi-stable model $f'
: X' \to C$, with $X'$ non-singular. By the definition of stable
surfaces in \cite{KS} (or of stable canonically polarized
varieties in \cite{KK}),
$$
f_* \omega^{[\eta]}_{X/C} = f'_* \omega^{\eta}_{X'/C}.
$$
\cite{V 2}, 3.2, gives an injective map
$$
f'_* \omega^{\eta}_{X'/C} \to \gamma^*  g_* \omega^{\eta}_{Y/B}
$$
which is an isomorphism over $\gamma^{-1}(B-S)$.
Hence
$$
{\rm deg} \ \Phi^* \lambda = \frac{{\rm deg} (f_*
\omega^{[\eta]}_{X/C} )^p}{{\rm deg} \ \gamma} \leq {\rm deg}
(g_* \omega^{\eta}_{Y/B})^p
$$
is bounded, as well.
\end{proof}

Let us denote by ${\rm \bf H} = {\rm \bf Hom} ((B, B-S),
(\bar{M}_h , M_h))$
the scheme parameterizing morphism
$\Phi : B \to \bar{M}_h$ with $ \Phi (B-S) \subset
M_h. $
Since $\lambda$ is ample, \ref{b_s} and \ref{b_cp} imply:
\begin{corollary} \label{b_maps}
Under the assumptions made in \ref{b_s} (or \ref{b_cp}) there
exists a subscheme $T \subset {\rm \bf H}$, of finite type over
$k$, which contains all points $[\Phi] \in {\rm \bf H}$, induced
by smooth morphisms
$$
g_0 : Y_0 \to B -S \in \sM_h (B - S).
$$
\end{corollary}

Without assuming the minimal model conjecture $MMP ({\rm deg} \
h + 1)$ \ref{b_cp} and \ref{b_maps} remain true, if $S =
\emptyset$. In fact, the existence of the very ample
sheaf $\lambda$ on $M_h$ has been shown in \cite{V 1}. For the
corresponding embedding $M_h \to \P^m$ choose $\bar{M}_h$ to
be the closure of $M_h$ in $\P^m$. Then for complete curves $B$ and
for morphisms $B \to \bar{M}_h$ with image in $M_h$, the arguments
used to prove \ref{b_cp} remain valid.

In \cite{V 1}, section 8.5, one finds the definition of a moduli functor
$\sD^{[N_0]}_{h}$ of canonically polarized normal varieties with
canonical singularities of index $N_0$. Kawamata \cite{Kaw} has
shown that this moduli functor is locally closed. Unfortunately
it is not known to be bounded, i.e. whether Matsusaka's big
theorem holds true.

As in \cite{V 1}, 1.20, one can enforce boundedness by
considering the submoduli functor $\sD^{[N_0] (M)}_{h}$ with
$$
\sD^{[N_0] (M)}_{h} (k) = \{ X \in \sD^{[N_0]}_{h} ;
\omega^{[N_0 \cdot M]}_{X} \ \ \mbox{very ample} \}.
$$
Then there exists a moduli scheme $D^{[N_0](M)}_{h}$ for
$\sD^{[N_0](M)}_{h}$. For $\nu = N_0 \cdot M$ and some $p > 0$
there exists again a very ample invertible sheaf $\lambda $ on
$D^{[N_0](M)}_{h}$ which induces
$$
{\rm det} (f_* \omega^{[\nu]}_{X/Z})^p
\mbox{ \ \ \ for \ \ \ } f: X \to Z \in \sD^{[N_0](M)}_{h} (Z).
$$
Choosing $\bar{M}_h$ again as the closure of $D^{[N_0](M)}_{h}$
for an embedding $\phi:D^{[N_0](M)}_{h}\to\P^m$
with $\phi^*\sO_{\P^m}(1)=\lambda$, one obtains:

\begin{corollary} \label{b_maps2}
For ${\rm deg} \ h > 2$, there exists a subscheme $T \subset
{\rm \bf H}$, of finite type over $k$, which contains all points
$[\Phi] \in {\rm \bf H}$, induced by morphisms
$$
g: Y \to B \in \sD^{[N_0](M)}_{h} (B),
$$
smooth over $B - S$.
\end{corollary}

\bibliographystyle{plain}

\begin{thebibliography}{XXX} 
\bibitem{Al} Alexeev, V.: Boundedness and $K^2$ for log
surfaces. Int. Journal Math. {\bf 5} (1994) 779 - 810.
\bibitem{A} Arakelov, A.: Families of algebraic curves with
fixed deneracies. Izv. Ak. Nauk. S.S.S.R., ser. Math. {\bf 35}
(1971) [Math. U.S.S.R. Izv. {\bf 5} (1971) 1277 - 1302].
\bibitem{EV 2} Esnault, H.; Viehweg, E.: Effective bounds for
semipositive sheaves and for the height of points on curves over
complex function fields. Compositio Math. {\bf 76} (1990) 69 -
85.
\bibitem{EV 1} Esnault, H.; Viehweg, E.: Lectures on Vanishing
Theorems. DMV Seminar {\bf 20} (1992), Birkh\"auser.
\bibitem{F} Faltings, G.: Arakelov's Theorem for abelian
varieties. Invent. math. {\bf 73} (1983) 337 - 348.
\bibitem{KK} Karu, K.: Minimal models and boundedness of stable
varieties. J. Alg. Geom. to appear.
\bibitem{Kaw} Kawamata, Y.: Deformations of canonical singularities, preprint
(1997).
\bibitem{Ko} Koll\'ar, J.: Projectivity of complete moduli. J.
Diff. Geom. {\bf 32} (1990) 235 - 268.
\bibitem{KS} Koll\'ar, J.; Shepherd-Barron, N. I.: Threefolds
and deformations of surface singularities. Invent. math. {\bf
91} (1988) 299 - 338.
\bibitem{K 2} Kov\'acs, S.: Smooth families over rational and
elliptic curves. J. Alg. Geom. {\bf 5} (1996) 369 - 385.
\bibitem{K 1} Kov\'acs, S.: On the minimal number of singular fibres
in a family of surfaces of general type. Journ. Reine Angew. Math.
{\bf 487} (1997) 171 - 177.
\bibitem{K 3} Kov\'acs, S.: Algebraic hyperbolicity of fine
moduli spaces. preprint (1998).
\bibitem{M 1} Migliorini, L.: A smooth family of minimal
surfaces of general type over a curve of genus at most one is
trivial. J. Alg. Geom. {\bf 4} (1995) 353 - 361.
\bibitem{V 2} Viehweg, E.: Weak positivity and the additivity of
the Kodaira dimension for certain fibre spaces. Adv. Stud. Pure
Math. {\bf 1} (1983) 329 - 353, North-Holland.
\bibitem{V 3} Viehweg, E.: Weak positivity and the additivity of
the Kodaira dimension II: The local Torelli map. Progr. Math.,
{\bf 39} (1983) 567 - 589, Birkh\"auser.
\bibitem{V 1} Viehweg, E.: Quasi-projective Moduli for Polarized
Manifolds. Ergebnisse 3. Folge, {\bf 30} (1995) Springer.
\bibitem{Z 1} Zhang, Qi: Holomorphic one-forms on projective
manifolds. J. Alg. Geom. {\bf 6} (1997) 777 - 787.

\end{thebibliography}

\end{document}